\documentclass[12pt]{article}
\newcommand\useplain[1]{#1}

\usepackage{amsthm}
\usepackage{fullpage}
\usepackage{amssymb, amsmath}
\usepackage{hyperref}
\usepackage{graphicx}

\newtheorem{theorem}{Theorem}
\newtheorem{thm}[theorem] {Theorem}
\newtheorem{lemma}[theorem]{Lemma}

\newtheorem{proposition}[theorem]{Proposition}

\newtheorem{problem}[theorem]{Problem}

\newtheorem{observation}[theorem]{Observation}
\theoremstyle{definition}

\newcommand\ex{\ensuremath{\mathrm{ex}}}
\theoremstyle{remark}

\useplain{\title{Multicolor Ramsey numbers for triple systems}}

\begin{document}
\author{
\quad{Maria Axenovich}
\thanks{Iowa State University, Ames,  IA, USA and
Karlsruher Institut f\"ur Technologie, Karlsruhe, Germany;  {\tt maria.aksenovich@kit.edu}, supported in part by NSF  grant DMS-0901008.}
\quad{Andr\'as Gy\'arf\'as}
\thanks{R\'enyi Institute of Mathematics, Budapest, Hungary;  {\tt gyarfas@renyi.hu}, supported in part by OTKA K104343.}\\
\quad{Hong Liu}
\thanks{
University of Illinois at Urbana-Champaign, Urbana, IL, USA; {\tt hliu36@illinois.edu}.}
\quad{Dhruv Mubayi}
\thanks{
University of Illinois at Chicago, Chicago IL, USA;  {\tt mubayi@uic.edu}, supported in part by NSF grant 0969092.}
}

\maketitle

\begin{abstract}
Given an $r$-uniform hypergraph $H$, the multicolor Ramsey number $r_k(H)$ is the minimum $n$ such that every $k$-coloring of the edges of the complete $r$-uniform hypergraph $K_n^r$ yields a monochromatic copy of $H$.
We investigate $r_k(H)$ when $k$ grows and $H$ is fixed.
For nontrivial  $3$-uniform hypergraphs $H$, the function $r_k(H)$ ranges from $\sqrt{6k}(1+o(1))$  to double exponential   in $k$.

We observe that $r_k(H)$ is polynomial in $k$ when $H$ is $r$-partite and at least single-exponential in $k$ otherwise. Erd\H{o}s, Hajnal and Rado gave bounds for large cliques $K_s^r$ with $s\ge s_0(r)$, showing its correct exponential tower growth. We give a proof for cliques of all sizes, $s>r$, using a slight modification of the celebrated stepping-up lemma of Erd\H{o}s and Hajnal.


For $3$-uniform hypergraphs, we give an infinite family with sub-double-exponential upper bound and
show connections between graph and hypergraph Ramsey numbers.  Specifically, we prove that
$$r_k(K_3)\le r_{4k}(K_4^3-e)\le r_{4k}(K_3)+1,$$ where $K_4^3-e$ is obtained from $K_4^3$  by deleting an edge.

We provide some other bounds, including single-exponential bounds for $F_5=\{abe,abd,cde\}$ as well as asymptotic or exact values of $r_k(H)$ when $H$ is the bow $\{abc,ade\}$, kite $\{abc,abd\}$, tight path $\{abc,bcd,cde\}$ or the windmill $\{abc,bde,cef,bce\}$. We also determine many new ``small'' Ramsey numbers and show their relations to designs. For example, the lower bound for $r_6(kite)=8$ is demonstrated by decomposing the triples of $[7]$ into six partial STS (two of them are Fano planes).

\end{abstract}

\parindent=0pt
\parskip=8pt

\section{Introduction, results} An $r$-uniform hypergraph $H$ is a pair $(V, E)$ where $V$ is a  vertex set and $E\subseteq \binom{V}{r}$ is the set of edges.  Let $K_n^{r}$ be the complete
$r$-uniform hypergraph containing all $r$-subsets of vertices as edges.
For an edge $\{v_1, v_2, \ldots, v_r\}$ we often write $v_1v_2\ldots v_r$.   When $r=2$,   denote $K_n^r$  by $K_n$. We shall also use the notation $\binom{[n]}{r}$ or $\binom{V}{r}$ for the edge set of  $K_n^r$.
An $r$-uniform hypergraph $H$ is $\ell$-partite if its vertex set can be partitioned into $\ell$ parts (called partite sets)  such that each edge contains at most  one
vertex from each part; $H$ is a complete $r$-partite hypergraph  if  each choice of
$r$ vertices from distinct partite sets forms an edge, and $H$ is balanced if its partite sets differ in size by at most one. A matching is a hypergraph consisting of disjoint edges.
A hypergraph $H=(V, E)$ is a subhypergraph of $F=(V', E')$ if $V\subseteq V'$ and $E\subseteq E'$. Denote by $\ex(n,H)$ the maximum number of edges in an $n$-vertex $r$-uniform hypergraph containing no copy of $H$ as subhypergraph. The density of an $r$-uniform hypergraph $H=(V,E)$ on $n$ vertices
is $d(H)=|E|/\binom{n}{r}$.

The {\bf multicolor Ramsey number} for an $r$-uniform hypergraph $H$, denoted by $r_k(H)$, is the minimum $n$ such that no matter how the edges of $K_n^r$ are colored with $k$ colors, there is a monochromatic copy of $H$.
While there is a number of results in the literature about $r_k(H)$ when $k$ is a small fixed number (see~\cite{C-F-S}), the case when $H$ is fixed and $k$ grows appears not to have been extensively studied.
The following three results are among the few results known in this area:

 \begin{theorem}[Lazebnik and Mubayi~\cite{L-M}]\label{Hst}
Fix integers $r,s,t\ge 2$. Let $H^r(s,t)$ be the complete $r$-partite $r$-uniform hypergraph with $r-2$ parts of size 1, one part of size $s$ and one part of size $t$. Then

(i) $r_k(H^r(2,t+1))=tk^2+O(k)$;

(ii) $r_k(H^r(s,t))=\Theta(k^s)$, for fixed $t,s\ge 2, t>(s-1)!$;

(iii) $r_k(H^r(3,3))=(1+o(1))k^3$.
\end{theorem}
Let $M$ be a matching with two $r$-tuples. Notice that an edge-coloring of $K_n^r$ without monochromatic copies of $M$ corresponds to a proper vertex-coloring of Kneser graph $K(n,r)$, that is, the graph with vertex set $\binom{[n]}{r}$ and two $r$-sets are adjacent if and only if they are disjoint.  Lov\'asz proved that the chromatic number of $K(n,r)$ is equal to $n- 2r+2$.  Reformulating his result, we obtain the following.
\begin{thm} [Lov\'asz~\cite{L}]\label{kneser}
If  $M$ is  a matching with two $r$-tuples, then $r_k(M)  =  k+2r -1$.
\end{thm}
Gy\'arf\'as and Raeisi observed that results of Cs\'ak\'any and Kahn~\cite{CSK} and the standard coloring of the Kneser graph imply the following.

\begin{proposition}[Gy\'arf\'as and  Raeisi~\cite{GYR}]\label{GR}
 If $C_3^3$ is the hypergraph with edge set $\{abc, cde, efa\}$, then $k+5 \leq r_k(C_3^3) \leq 3k+1$.
 \end{proposition}

In this paper,  we start  a systematic investigation on the growth rate of $r_k(H)$ for some fixed $H$ as $k$ grows.
Our first result shows that $r_k(H)$ is polynomial in $k$ if and only if $H$ is $r$-partite.

\begin{proposition}\label{non-r-partite}
Let $r\ge 2$ be  fixed  and  $H$ be a connected $r$-uniform hypergraph.
Then $r_k(H)$ is polynomial in $k$ if and only if $H$ is $r$-partite. In particular, there are positive constants $c$ and $c'$, such that

(i) If $H$ is $r$-partite, then $r_k(H)=O(k^{c})$

(ii) If $H$ is not $r$-partite, then $r_k(H)\ge 2^{c'k}.$
\end{proposition}

Determining the growth rate of $r_k(H)$ in general is known to be a very hard problem. For example, the best known bounds even for the smallest nontrivial graph case are $c^k<r_k(K_3)<c'k!$ for some positive constants $c$ and $c'$ (see Chung~\cite{Chung} and Erd\H{o}s, Szekeres~\cite{E-Sz}). Define the tower function as follows: $t_1(n)=n$ and $t_{i+1}(n)=2^{t_{i}(n)}$ for all $i\ge 1$. Erd\H{o}s, Hajnal and Rado gave an upper bound for all cliques and a lower bound for only large cliques.

\begin{theorem}[Erd\H{o}s and Rado~\cite{E-R}, Erd\H{o}s et al.~\cite{E-H-R}]\label{clique-upperbound}
Let $s>r\ge 2$.  There are positive integers $c=c(s,r)\le 3(s-r)$, $s_0(r)$,  and $c'=c'(s,r)$ such that
$$t_r(c'k)< r_k(K_s^{r})<  t_{r}(ck\log k)$$
where the lower bound holds for $s\ge s_0(r)$.
\end{theorem}

It is worth noting that the lower bound in~\cite{E-H-R} was stated for the case when the number of colors, $k$, is fixed while $r$ grows and the bound was only for large cliques. But the proof in~\cite{E-H-R} applies naturally to our case as well, when $k$ grows and the other parameters are fixed. Recently, an improved stepping-up lemma was proved by Conlon et al~\cite{C-F-S-steppingup}. Their main result implies a lower bound for cliques of smaller sizes, but still only for $s\ge 2r-1$. Duffus, Lefmann and R\"odl~\cite{DLR} took another approach, using shift graphs, and proved a lower bound for cliques of all sizes $s> r$, but require $k$ being fixed and $r\gg k$. Our next result gives a proof for cliques of all sizes using a slight modification of the stepping-up lemma, due to Erd\H{o}s and Hajnal (see Chapter 4.7 in~\cite{G-R-S}).

\begin{theorem}\label{thm-complete}
For any $s>r\ge 2$ and $k>r2^r$ we have
$$r_k(K_s^{r})> t_{r}\left(\frac{k}{2^r}\right).$$
\end{theorem}

Our remaining results are all for $3$-uniform hypergraphs and we will address the question of determining $r_k(H)$ for most interesting $H$'s with $6$ or fewer vertices. Let $K_4^3-e$ be a hypergraph obtained from $K_4^3$ by removing one edge. Our next theorem gives bounds on $r_k(K_4^3-e)$ in terms of $r_k(K_3)$, showing that compared to the double-exponential bounds for $K_4^3$ from Theorems~\ref{clique-upperbound} and~\ref{thm-complete}, the correct order of magnitude for $r_k(K_4^3-e)$ is single-exponential.

\begin{theorem}\label{K43-} For any $k \geq 2$,
$$r_k(K_3)\le r_{4k}(K_4^{3}-e)\quad\mbox{  and  }\quad r_{k}(K_4^{3}-e) \le   r_{k}(K_3)+1.$$
Moreover $r_2(K_4^3-e) = r_2(K_3)+1 =7$.
\end{theorem}

Denote by $F_5$ the hypergraph with edges $\{abc,abd,cde\}$.
We show that $r_k(F_5)$ behaves similarly to  $r_k(K_3)$.

\begin{theorem}\label{F5}
There is a  positive constant $c$ such that, for $k\ge 4$,
$$2^{ck} \leq r_k( F_5)\leq k!,$$
and $r_2(F_5)=6, r_3(F_5)=7$.
\end{theorem}

The simplest non-trivial  triple systems have just two edges. The {\bf kite} is a $3$-uniform hypergraph  with two edges sharing two vertices.
The {\bf  bow} is a $3$-uniform hypergraph with  two edges sharing a single vertex.

\begin{theorem}\label{bow} Let  $r_k=r_k(bow)$.  Then
$$r_k= (1+o(1))\sqrt{6k}.$$
If $k={{n\choose 3}\over n}$ and  $n\equiv 4,8 \pmod{12}$, then $r_k=n+1$. Moreover, $r_2=5,r_3=r_4=r_5=6$, $r_6=7$, $r_7=r_8=r_9=r_{10}=9$,
$9\le r_{11}\le r_{12}\le r_{13}\le r_{14}\le 10, r_{15}=11.$
\end{theorem}

{\bf Remark.}
Note that $r_k(bow)$ is the smallest multicolor Ramsey number among nontrivial  $3$-uniform hypergraphs
since $r_k(H) \geq \min\{ r_k(bow),  r_k(kite), r_k(M)\}$, where
$M$ is a matching with $2$ triples. Indeed, each nontrivial  $3$ uniform hypergraph contains at least two edges that
form one of $bow, kite$ or $M$, and  Theorem~\ref{kneser} gives $r_k(M)=k+5$.

\begin{theorem}\label{kite} Let $r_k=r_k(kite)$. Then
$$r_k=
\begin{cases}
k+1, & \mbox{if } k\equiv 3\mbox{ (mod 6)} \\
k+1~\mbox{ or }~k+2, & \mbox{if } k\equiv 4\mbox{ (mod 6)} \\
k+2, & \mbox{if } k\equiv 0,2\mbox{ (mod 6)} \\
k+3, & \mbox{if } k\equiv 1,5\mbox{ (mod 6)},  k\ne 5 \\
6 & \mbox{if } k=5,\\
5 & \mbox{if } k=4 \end{cases}$$
\end{theorem}

Let $a,b$ be positive integers. Denote by $F(a,b)$ the $3$-uniform hypergraph with vertex set $V=A\cup B$,  $A\cap B = \emptyset$, $|A|=a, |B|=b$ and edge set consisting of all triples with one vertex in $A$ and two vertices in $B$ (for example, $F(2,2)$ is the kite).

\begin{proposition}\label{fa2}
For any $a\ge 2$, we have
$$k(a-1) < r_k(F(a,2))\le k(a-1)+3.$$
\end{proposition}

In general, $r_k(F(a,b))$ grows slower than double exponential in $k$ and possibly faster than exponential in $k$. (Recall that Theorems~\ref{clique-upperbound} and~\ref{thm-complete} give double-exponential bounds.)

\begin{theorem}\label{fab}
Given $3\le a\le b$, we have, for positive constants $c=c(a,b)$ and $c'=c(a,b)$
$$2^{c'k} <r_k(F(a,b))<r_t(K_b)+m< 2^{ck^{a+1}\log k} ,$$
where $m=(a-1)k+1$, and $t=k{m\choose a}$.
\end{theorem}

The {\bf windmill} $W$ with \emph{center edge} $abc$ is the hypergraph with six vertices and edge set $\{abc,abd,bce,acf\}$.

\begin{thm}\label{windmill}
$$(1-o(1))3k\le r_k(W)\le 3k+3.$$
\end{thm}

It is interesting to compare Theorem \ref{windmill} with Proposition \ref{GR}.
In fact, the upper bounds in both cases come  from the corresponding Tur\'an-type results. Indeed, $\ex(n,C_3^3)={n-1\choose 2}$  (Frankl-F\"uredi \cite{FF} for large $n$, Cs\'ak\'any-Kahn~\cite{CSK} for $n\ge 6$) while  $\ex(n,W)\le {n\choose 2}$  (\cite{FF}).

The ideas giving the asymptotic of $r_k(W)$ can be also used for the {\bf tight path } $P_3^3=\{abc,bcd,cde\}$.

\begin{theorem}\label{tightp}
$2k(1-o(1))\le r_k(P_3^3) \leq 2k+3.$
\end{theorem}

The rest of the paper will be organized as follows. In Section~\ref{sec-gb}, we give some auxiliary results and prove Proposition~\ref{non-r-partite}. Theorems ~\ref{thm-complete} - \ref{tightp} will be proved in Sections~\ref{sec-complete}-\ref{windpasch}.  Section ~\ref{small} is devoted to exact values of Ramsey numbers for small number of colors and Section ~\ref{conclude} contains remarks, conjectures and problems.

In some later sections we give lower bounds on Ramsey numbers based on block designs.  A $t-(v,k,\lambda)$ {\bf design} is a  subset of $\binom{[v]}{k}$, called blocks,  such that each $t$ element subset of $[v]$ is
contained in exactly $\lambda$ blocks.

\section{General bounds and auxiliary results}\label{sec-gb}
In this section we prove some general bounds on $r_k(H)$ and obtain some consequences including Proposition \ref{non-r-partite}.
Recall that  the density of an $r$-uniform hypergraph $F$ with $n$ vertices and $e$ edges is $d(F)= \frac{e}{\binom{n}{r}}$.

\begin{lemma}\label{general}
Let $H$  be a fixed $r$-uniform hypergraph and $F$ be  an $r$-uniform hypergraphs with $n$ vertices,  density $d(F)=d$, and  not containing copies of $H$ as a
subhypergraph. Then

(i) $r_k(H) \leq  1 + \max\{ n:   \lceil \binom{n}{r} / \ex(n, H) \rceil \leq k \}$,\\
(ii) If  $\binom{n}{r}(1-d)^k < 1$ then $r_k(H) \geq n$.
\end{lemma}

 {\bf Proof.}
(i)  Consider a coloring of  $K_n^r$ with $k$ colors and no monochromatic copy of $H$.
Then each color class has at most  $\ex(n,H)$ edges.

(ii) Consider $k$ copies of  hypergraph $F$ obtained by mapping its vertices randomly to a given  set $V$ of $n$ vertices.
Here, we choose vertex permutations uniformly.  Assign  the edges of the $i$th copy  of $F$ color $i$, $i=1, \ldots, k$.
If an edge belongs to several copies of $F$, assign the smallest available label.
We claim that with positive probability, each edge of $K= \binom{V}{r}$ belongs to some copy of $F$.
Indeed, the probability that a given  edge  of $K$  uncovered is $(1-d)^k$.
Thus, the probability that there is an uncovered edge of $K$ is at most $\binom{n}{r}(1-d)^k<1 $.
Therefore, with positive probability, all edges are covered and the resulting coloring of $K$   contains no monochromatic copy of $H$.
\qed
\medskip

 {\bf Proof of Proposition~\ref{non-r-partite}.}  (i) The proposition follows from Lemma~\ref{general}(i) by using the fact that
 $\ex(n,H)<n^{r-c}$ for some positive constant $c=c(H)$, when $H$ is $r$-partite, see \cite{E}.
 So,  $k \geq \binom{n}{r}/ \ex(n, F) \geq Cn^r/n^{r-c} = Cn^c$, for a constant $C=C(r)$.
Thus $n\leq C^{-1/c} k^{1/c}$.

(ii) Let $H$ be non-$r$-partite.   Apply Lemma~\ref{general}(ii) with   $F$ being    a complete $r$-uniform $r$-partite balanced hypergraph   on $n=2^{c'k}$ vertices (and $r|n$).
Clearly $H$ is not contained in $F$ as a subgraph.  Moreover,
$ d(F)\ge \frac{(n/r)^r}{\binom{n}{r}}>\frac{(n/r)^r}{(en/r)^r}=e^{-r}.$
Hence  for $k =  c\log n$ and $c>e^r(r+1)$,
$$\binom{n}{r}(1-d)^k =  \binom{n}{r}(1-d)^{c\log n} <n^re^{-cd\log n}=e^{(r-cd)\log n}<1.\quad \quad
\qed $$

The \emph{trace} of a $3$-uniform hypergraph $H$ at vertex $v$ is the graph  on vertex set $V(H)-\{v\}$ and with edge set $\{e-\{v\}:~  e\in H,  ~ v\in e\}$. A transversal of a hypergraph is a set of vertices non-trivially intersecting each edge.

\begin{lemma}\label{single-transversal}
Let $H$ be a 3-uniform hypergraph with a single-vertex  transversal $\{v\}$.
Let $G$ be  a trace of $H$ with respect to $v$.
Then $r_k(H) \leq r_k(G)+1$.
\end{lemma}

 {\bf Proof.}
Given a $k$-coloring $c$ of ${[n] \choose 3}$ with no monochromatic $H$, let $c'$ be the $k$-coloring of ${[n-1] \choose 2}$ defined by $c'(ij)=c(ijn)$.  Then $c'$ has no monochromatic $G$ and consequently $r_k(G) \ge r_k(H)-1$ as required. \qed

\section{$K_s^{r}$ for $s>r\ge 2$}\label{sec-complete}
In this
section we prove Theorem \ref{thm-complete} using a variant of the stepping-up lemma of Erd\H os and Hajnal.


\noindent {\bf Proof of Theorem \ref{thm-complete}.}
It suffices to prove the result for $s=r+1$ since $r_k(K_s^{r})\ge r_k(K_{r+1}^{r})$ for any $s>r$.
We use induction on $r$ to show that  $r_k(K_{r+1}^r)\ge t_r(k/2^{r-2} - 2r)$ for all $k \ge r2^r$. Since $k\ge r2^r$, we have $k/2^{r-2}-2r \ge k/2^r$ and the result follows.

The base case $r=2$ is given by $r_k(K_3)> 2^k>2^{k-4}=t_2(k-4)$. Assume the result holds for some $r\ge 2$ and let $n=r_k(K_{r+1}^{r})-1$.  By the inductive hypothesis
$$n\ge t_{r}(k/2^{r-2}-2r)-1.$$ Let $\phi: \binom{[n]}{r}\rightarrow [k]$ be a coloring with no monochromatic $K_{r+1}^{r}$.
We will construct a coloring $\psi: \binom{[2^n]}{r+1}\rightarrow [2k+2r-4]$ with no monochromatic  $K_{r+2}^{r+1}$.
This shows that
$$r_{2k+2r -4}(K_{r+2}^{r+1})\ge 1+2^n\ge 1+
 \frac{1}{2}t_{r+1}(k/2^{r-2}-2r).$$  Now suppose we are given $k'\ge (r+1)2^{r+1}$.  If $k'-2r+4$ is odd, then let $k''=k'-1$
and if $k'-2r+4$ is even then let $k''=k'$. Set $k=(k''-2r+4)/2$ (which is an integer) and observe that
$k\ge r2^r$ and $k''=2k+2r-4$.
 Then
 $$k/2^{r-2}-2r\ge k''/2^{r-1} - 2(r+1) +1$$
 and $r_{k'}(K_{r+2}^{r+1})$ is at least
 $$r_{k''}(K_{r+2}^{r+1})\ge 1+
 \frac{1}{2}t_{r+1}(k/2^{r-2}-2r)\ge 1+\frac{1}{2}
 t_{r+1}(k''/2^{r-1} - 2(r+1) +1)>t_{r+1}(k'/2^{r-1} - 2(r+1)).$$
 Now we shall construct a coloring $\psi$ of  $\binom{[2^n]}{r+1}$ using the coloring $\phi$ of  $\binom{[n]}{r}$ that has no monochromatic $K_{r+1}^{r}$.
Represent the elements of $[2^n]$  with  $0$-$1$-sequences on $n$ coordinates. For a vertex $u$ and integer $i$, we denote $u(i)$
the $i$th coordinate of $u$ in this representation.
 Given two vertices $u,v\in [2^n]$,  say that $u<v$  if $u(i)<v(i)$ and $u(j)=v(j)$ for $j<i$. Denote such an $i$ by $f(uv)$. Given any $u_1<\cdots <u_{r+1}$, let $f_i:=f(u_iu_{i+1})$, for every $1\le i\le r$. Observe crucially that

\medskip
(1) $f_i\neq f_{i+1}$, for every $1\le i\le r-1$;

\medskip

(2) $f(u_1u_{r+1})=\displaystyle\min_{1\le i\le r}\{f_i\}$ and the minimum is reached by a unique $i$.

\medskip

We define coloring $\psi$ as follows:
 \begin{equation}\nonumber
\psi(u_1...u_{r+1})= \begin{cases}
  (\phi(f_1,..., f_{r}), 1) &  \mbox{ if $(f_1,...,f_r)$ is an increasing sequence}, \\
  (\phi(f_1,..., f_{r}), 2) &  \mbox{ if $(f_1,...,f_r)$ is a  decreasing sequence},\\
  (i,3) &  \mbox{ if $f_1<f_2<\cdots <f_i>f_{i+1}$}, 2\le i\le r-1, \mbox{ for }r\ge 3, \\
  (i,4) & \mbox{ if $f_1>f_2>\cdots >f_i<f_{i+1}$}, 2\le i\le r-1, \mbox{ for }r\ge 3.
 \end{cases}
 \end{equation}

Suppose to the contrary that there is a monochromatic copy of $K_{r+2}^{r+1}$ under $\psi$ on vertex set $U=\{u_1,...,u_{r+2}\}$ with $u_1<\cdots<u_{r+2}$. Without loss of generality, we distinguish two cases.

\medskip

\noindent\textbf{Case 1:} The second coordinate of $\psi$ on each $(r+1)$-tuple is $1$. First notice that the second coordinate of $\psi$ on $u_1,...,u_{r+1}$ and $u_2,...,u_{r+2}$ being $1$ implies $f_1<f_2<\cdots <f_{r}<f_{r+1}$ and together with (2), we have $f(u_1u_i)=f(u_1u_2)=f_1$ for all $3\le i\le r+2$. Similarly from $u_2,...,u_{r+2}$, we have that for every $2 \le p<q\le r+2$, $f(u_pu_q)=f_p$.  Recall that the color of the $(r+1)$-set $\{u_1,...,u_{r+2}\}- \{u_i\}$ under $\psi$ is determined by the color of the $r$-set  $\{f_1,...,f_{r+1}\}-\{f_i\}$ under $\phi$.
 Let $F:=\{f_1,...,f_{r+1}\}$ and $U=\{u_1, \ldots, u_{r+2}\}$.
 Let us denote the above implication by  $$U\setminus\{u_i\}\Rightarrow F\setminus\{f_i\}.$$
Thus a monochromatic $K_{r+2}^{r+1}$ on $U$ under $\psi$ yields a monochromatics $K_{r+1}^{r}$ on $F$ under $\phi$, a contradiction.

 \medskip

\noindent\textbf{Case 2:} Each $(r+1)$-tuple gets color $(i,3)$  for some $i$ with $2\le i\le r-1$. Then  $\psi(u_1,...,u_{r+1})=(i,3)$ implies $f_i>f_{i+1}$. On the other hand, $\psi(u_2,...,u_{r+2})=(i,3)$ implies $f_i<f_{i+1}$, a contradiction.

If the second coordinate is 2 or 4 the arguments are almost identical  to those in  Case 1 or 2.
 \qed

\section{$K_4^{3}-e$ and $F_5$}\label{sec-K43}

Notice that in contrast to the double-exponential growth for $K_4^3$, $r_k(K_4^3-e)$ is single-exponential in the number of colors $k$. Indeed, since $K_4^3-e$ is not $3$-partite, Proposition~\ref{non-r-partite} yields $r_k(K_4^3-e)>2^{ck}$. For the upper bound, one can use a variation of the classical Erd\H{o}s-Rado pigeonhole argument to obtain $r_k(K_4^3-e)<2^{(k+1)\log k}$. We will, however, use a different approach to prove this fact, which also shows some connection between the multicolor Ramsey number of $K_4^3-e$ and the multicolor Ramsey number of a triangle.

 {\bf Proof of Theorem \ref{K43-}.}
For the lower bound, let $n=r_k(K_3)-1$ and $\phi: \binom{[n]}{2}\rightarrow k$ be a $k$-coloring of $\binom{[n]}{2}$ with no monochromatic triangles. We will construct a coloring $\psi$ of
 $\binom{[n]}{3}$ with $4k$ colors with  no monochromatic  $K_4^3-e$. This then would imply that  $r_{4k}(K_4^3-e)\ge n+1=r_k(K_3)$ as desired. Let $\psi$ be  the following coloring of the triples $i<j<k$.  If $P$ is a path with vertices $i,j,k$, denote by $\phi'(P)$ the color under $\phi$ of the edge in $\{i,j,k\}$ that is not in $P$.
For such a path $P$, let the type of $P$, $t(P)=1, 2$, or $3$  if $i, j$ or $k$  is its center, respectively.  If $\{i,j,k\}$  is a rainbow triangle, let  $\psi(ijk)=(0,\phi(jk))$.  If $ \{i,j,k\}$  induces  a monochromatic  path $P$, let $\psi(ijk)=(t(P), \phi'(P))$.

Suppose there is a monochromatic copy $K=\{abc, abd, acd\}$ of $K_4^3-e$, we will show a contradiction when the first coordinate is $0$, namely all three triples $\{abc,abd,acd\}$ span rainbow triangles under $\phi$. The cases when the first coordinate is $1,2$ or $3$, can be proved using a similar argument. Notice that when the first coordinate is $0$, by the definition of $\psi$, the color of a triple depends on the color, under $\phi$, of the edge spanned by the two largest elements in that triple.
Since $b,c,d$ play a symmetric role, we can assume that $b<c<d$. If $a$ is the smallest, then $\psi(abc)=\psi(abd)=\psi(acd)$ implies $\phi(bc)=\phi(bd)=\phi(cd)$, i.e. $bcd$ is monochromatic under $\phi$. Thus $b$ is the smallest. But then $\psi(abc)=\psi(abd)$ implies $\phi(ac)=\phi(ad)$, which means $acd$ is not a rainbow triangle under $\phi$, a contradiction.

%
%
%
%
For the upper bound, simply notice that $K_4^3-e=\{abc,abd,acd\}$ has a single vertex transversal $\{a\}$, and the trace of $a$ is a triangle on $\{b,c,d\}$. Thus the upper bound follows from Lemma~\ref{single-transversal}.
The case with $2$ colors is treated in Section \ref{small}.
\qed

\medskip

\noindent {\bf Proof of Theorem \ref{F5}.} The cases $k=2,3$ are treated in Section \ref{small}. The general lower bound follows from Proposition~\ref{non-r-partite}(ii), since $F_5$ is not $3$-partite.

The upper bound follows by induction with basis $k=4$.
Suppose that the edges of $K_{24}^3$  with vertex set $V$ can be $4$-colored so that there is no monochromatic $F_5$. There are $22$  triples $uvx$ containing a
fixed pair $uv$.  Assume that  $uvx_1,uvx_2$ are
red triangles. Then $x_1x_2y$ cannot be red for $y \in Y=V-\{u,v,x_1,x_2\}.$
Thus we have a set $S$, $S\subseteq Y$,   $|S|\geq \lceil(|V|-4)/3 \rceil  = 7$ and
$x_1x_2y$ are blue triples for all $y\in S$.  Therefore, no triple in $S$ is colored blue, and thus $\binom{S}{3}$ uses
$k-1= 3$ colors.   Applying Theorem \ref{r3f5} to the $3$-colored subhypergraph spanned by $S$, we get a contradiction.

The inductive step is simply repeating the argument above in general. Suppose we already know $r_k(F_5)\le k!$ for some $k\ge 4$ and we have a $K_n^3$ with a $(k+1)$-coloring such that there is no monochromatic $F_5$. Selecting $u,v,x_1,x_2$ as above and applying the same argument, we get  $n-4\le k(k!-1) < (k+1)!-k$, thus $n\leq (k+1)! - k+4 \leq (k+1)!.$
This implies $r_{k+1}(F_5)\le (k+1)!$.  \qed

\medskip
{\bf Remark.}
The above results slightly suggests  that $r_k(F_5)\le r_k(K_3)$ might hold.
Although the bound $r_k(F_5)\le k!$ in Theorem \ref{F5} can be improved slightly, this improvement
still does not show that  $r_k(F_5)\le r_k(K_3)$.

\section{Bow,  Kite, $F(a,b)$}\label{sec-Bow}

The next lemma (without the statements on the extremal configurations) is referred in \cite{L-R} as an unpublished remark of Erd\H{o}s and S\'os.

\begin{lemma}\label{tur} $$\ex(n,bow)=\left\{\begin{array}{ll} n &\mbox{if $n\equiv 0 \pmod{4}$}\\
                    n-1 &\mbox{if $n\equiv 1 \pmod{4}$}\\
                    n-2 &\mbox{if $n\equiv 2,3 \pmod{4}.$}
\end{array}\right.$$  When $n\equiv 0,1 \pmod{4}$, the extremal configurations are unique, all components are $K_4^3$-s, (apart from a possible one vertex component). When $n\equiv 2 \pmod{4}$, the extremal configuration is either ${n-2\over 4}$ copies $K_4^3$-s and two isolated vertices or any number of $K_4^3$-s and one star component. Similarly, when $n\equiv 3 \pmod{4}$, the extremal configuration is either ${n-3\over 4}$ copies $K_4^3$-s and component with a single edge or any number of $K_4^3$-s and one star component.
\end{lemma}

\noindent {\bf Proof.}
Suppose $C$ is the vertex set of a nontrivial connected component of a $3$-uniform hypergraph without a $bow$. Then either $C$ spans only one edge or there are two edges $e_1,e_2$  in $C$, intersecting in two vertices, $u,v$. Suppose that $|C|>4$. Then every edge $f$ that is not covered by $e_1\cup e_2$ and intersecting $e_1\cup e_2$ must contain $u,v$ and a vertex $w$ not covered by $e_1\cup e_2$. It is easy to see that these  vertices $w$ cover $C$ and $C$ has no other edges, thus $C$ has $|C|-2$ edges, all containing $u,v$. Such a component is called a star component.

On the other hand, if $|C|=4$ then we have two, three or four edges in $C$. From this analysis the lemma follows.
\qed

Lower bounds of $r_k(bow)$ follow from the existence of resolvable designs. A $3-(n,4,1)$ design is a set of $4$-element subsets (blocks) of an $n$-element set $V$ such that each $3$-element subset of $V$ is in precisely one block. Hanani \cite{HA} showed that $3-(n,4,1)$ designs exist if and only if $n\equiv 2,4 \pmod{6}$. A $3-(n,4,1)$ design is called {\em resolvable} if its blocks can be grouped so that each group (parallel class) gives a partition of $V$.
Resolvable $3-(n,4,1)$ designs exist if and only if $n\equiv 4,8 \pmod{12}$, see \cite{HAR1,HAR}, and \cite{JZHU}.

\noindent {\bf Proof of Theorem \ref{bow}.}
When $n\equiv 4,8 \pmod{12}$,  $k={{n\choose 3}\over n}$,  $\ex(n,bow)=n$, thus  Lemma~\ref{general} (i) gives $r_k\le n+1$.
This is sharp, since $K_n^3$ can be partitioned into $k$ matchings. The statement $r_k(bow) \approx \sqrt{6k}$ follows from
considering  the construction for largest $n$, $n\equiv 4,8 \pmod{12}$,  $k \geq {{n\choose 3}\over n}$ for the lower bound and
applying the Lemma~\ref{general}(i) for the upper bound.  The statements about the small values are proved in Section \ref{small}. \qed

\medskip

\noindent {\bf Proof of Theorem \ref{kite}.}
Let $H=F(2,2)$ be the kite. Then  $\ex(n,H)$ corresponds to the maximum number of triples on $n$ elements such that any two triples intersect in at most one element, i.e. the maximum number of edges in a linear $3$-uniform hypergraph.
A well-known result of Sch\"onheim \cite{SCH} and others (the cases $n\equiv 0,1,2,3 \mbox{ (mod 6)}$ go back even to Kirkman \cite{KI}) is
$\ex(n, H)=\left\lfloor {n\over 3}\left\lfloor{n-1\over 2}\right\rfloor\right\rfloor - \epsilon,$
where  $\epsilon=1$ for $n\equiv 5 \mbox{ (mod 6)}$, otherwise $\epsilon=0$.
Lemma~\ref{general}(i)  gives, after some   calculations,  the upper bounds.

The lower bound for the cases $k\equiv 3,4\mbox{ (mod 6)}$ is easy.
Given $K_n^{3}=(V,E)$, consider $V=Z_n$ and color triple $ijk$ with color $i+j+k$ (mod $n$). Clearly this coloring yields no monochromatic $H$, hence $r_k(H)>k$.

For the cases $k\equiv 0,1,2,5\mbox{ (mod 6)}$ the (difficult) constructions of J. X. Lu \cite{LU1,LU2} finished by Teirlinck \cite{TEI} are needed: for $n>7,n \equiv 1,3\mbox{ (mod 6)}$, $K_n^3$ can be partitioned into $n-2$ Steiner triple systems (called a {\em large set} of STS).

Indeed, for $k\equiv 0,2\mbox{ (mod 6)}$ we need a kite-free $k$-coloring of $K_{k+1}^3$ i.e. $(n-1)$-coloring of $K_n^3$ when $n\equiv 1,3\mbox{ (mod 6)}$. This can be done even with $n-2$ colors according to the cited result of Lu. However, the case $k=6$ is exceptional because Lu's theorem does not hold for $n=7$. Nevertheless, there is a $6$-coloring of $K_7^3$ without a monochromatic kite as shown in Proposition \ref{r5kite}. Similarly, for $k\equiv 1,5\mbox{ (mod 6)}$ we need a kite-free $k$-coloring of $K_{k+2}^3$ i.e. $(n-2)$-coloring of $K_n^3$ when $n\equiv 3,1\mbox{ (mod 6)}$. This is provided by Lu's theorem, apart from the case $k\equiv 5$ ($n=7$) which is indeed exceptional, in Proposition \ref{r5kite} we prove that $r_5(kite)=6$ (together with the case $k=4$).
\qed

\medskip

\noindent {\bf Proof of Proposition \ref{fa2}.} In an $F(a,2)$-free coloring of $K_n^3$ any pair of vertices is in at most $a-1$ edges of the same color. Thus $n\le 2+k(a-1)$, proving the upper bound. (One can also use Lemma~\ref{single-transversal} and the multicolor Ramsey number for stars (see~\cite{B}): $r_k(K_{1,a})\le k(a-1)+2$.)

For the lower bound, set $n=k(a-1)$ and consider $K_n^{3}=(V,E)$ with $V=Z_n$. Color a each edge with the sum of its vertices mod $k$. Then a monochromatic copy of $F(a,2)$ would require that for some $y,z\in V$, $y+z+x_1,...,y+z+x_a$ are all equal (mod $k$) i.e. we have $a$ different positive $x_s$, all equal (mod $k$), which is impossible. Hence $r_k(F(a,2))>k(a-1)$.
\qed

\medskip

\noindent {\bf Proof of Theorem \ref{fab}.}
For the upper bound, let $N=r_t(K_b)+m$. Consider a $k$-coloring $\phi$ of the triples of $K_N$. Fix a set $S$ of $m$ vertices and define a $t$-coloring $c$ on the pairs of the remaining $N-m$ vertices as follows. Let $c(xy)=(\phi(xys_i),s_1,s_2,...,s_a)$, where $\phi(xys_i)$ is the majority color on triples containing $x$ and $y$, and $s_1,s_2,...,s_a
\in S$ is the lexicographically first $a$-tuple in $S$ such that $\phi(xys_i)=\phi(xys_j)$ for every $1\le j\le a$ (by the choice of $m$ there is
such an $a$-tuple). Since $c$ is a $t$-coloring of a complete graph on $N-m=r_t(K_b)$ vertices, there is monochromatic $K_b$ in $c$, which gives a monochromatic $F(a,b)$ in $\phi$.

 A lower bound for $r_k(F(a,b))$ is obtained from Proposition \ref{non-r-partite} (i) since $F(a,b)$ is not $3$-partite, for  $b\geq 3$.
\qed


\section{Windmill and tight path}\label{windpasch}

The following result (conjectured by Kalai) is a special case of a theorem of F\"uredi and Frankl (\cite{FF}, Theorem 3.8). We give their proof also, since it is extremely short in this special case.

\begin{thm}\label{exW}
$\ex(n, W)\le \binom{n}{2}$ with equality for every $n\equiv 1,5\mbox{ (mod 20)}$.
\end{thm}
\noindent {\bf Proof.}
The lower bound comes from the following construction. Let $n\equiv 1,5\mbox{ (mod 20)}$ and consider a Steiner system $S$, a $2-(n,5,1)$ design, i.e., a set of $5$-element blocks on $n$ elements such that every pair lies in precisely one block.  Its existence is proved by Hanani \cite{HA1,HA2}.
Then the number of blocks is $\binom{n}{2}/10$.
Now place $10$ triples inside  each block  of $S$. The resulting triple system, $H$, has ${n \choose 2}$ triples and  is $W$-free.
Indeed, a copy of $W$ would have to be contained in one of the blocks, but each block has less vertices than the number of vertices in $W$.

To prove the upper bound, suppose that $H$ is a $3$-uniform hypergraph with no $W$. For $x,y\in V(H)$, the {\em codegree}  $d(x,y)$ is the number of edges of $H$ containing both $x,y$.
Let   $a,b,c$ be codegrees of  three pairs  of vertices from a   edge of $H$, $1\leq a \leq b\leq c$.
If  $a=2$, $b\geq 3$ and $c\geq 4$, then $H$ contains a copy of  $W$.
Thus either $a=1$ or $a=b=2$ or $a=2, b=3, c=3$.  In each of these cases we have that
$1/a+1/b+1/c \geq 1$.  For each edge $e=uvw$ of $H$,  let
$$w(e) = \frac{1}{d(u,v)} +\frac{1}{d(v,w)}+\frac{1}{d(u,w)}.$$
We see that $w(e)\geq 1$.
 Let $s = \sum_{e\in H} w(e)$. Notice that $s\le {n \choose 2}$, since a term $\frac{1}{d(u,v)}$  appears exactly $d(u,v)$  times  for each pair $uv$ that belongs to at least one  edge of $H$.
 Now,  $|H|\le |H|\min_{e\in H} w(e)\le s \le {n \choose 2}.$
\qed

For the next proof we need the following decomposition result:
\begin{thm}[Pippenger and Spencer  ~\cite{P-S}]\label{pippsp}
  Let $r$ be fixed and $D$ be sufficiently large. Let $H$ be an $r$-uniform hypergraph with $d(v)=(1+o(1))D$ for every $v\in V(H)$ and codegree of
  $o(D)$ for every pair $\{u,v\} \subseteq V(H)$.  Then $E(H)$ can be partitioned into $(1+o(1))D$ matchings.
\end{thm}

\noindent {\bf Proof of Theorem \ref{windmill}.} To prove the lower bound, let $S$ be a $3-(n,5,1)$ design, i.e. a set of $5$-element blocks of an
$n$-element set such that each  $3$-element set is in precisely one block. The existence of such designs are known for infinitely many $n$, for example for $n=4^s+1, s\ge 2$ \cite{HU}, see also \cite{MRC}.
Construct an auxiliary $10$-uniform hypergraph $H$ where $V(H)$ is the set of ${n\choose 2}$ pairs in $V(S)$, and
ten of these pairs form an edge of $H$ if and only if they are the ten pairs in a block of $S$. Since every pair in $V(S)$ is in exactly $(n-2)/3$ blocks of $S$, $H$ is an $(n-2)/3$-regular hypergraph. On the other hand, the codegree of any two vertices in $H$ is at most one. Indeed, any two vertices in $H$ (two pairs in $V(S)$) contain at least three vertices in $V(S)$, and they can be in at most one block of $S$. With large enough $n$, and with $r=10, D=n/3$, the conditions of Theorem \ref{pippsp} hold so we can decompose $E(H)$ into $m=(1+o(1))n/3$ matchings $M_i, i=1,2,\dots,m$. Each $M_i$ corresponds to a subset of blocks $S_i$ of $S$ and any two blocks in $S_i$ share at most one element in $V(S)$.  The set of triples covered by the blocks of any $S_i$ form a $W$-free triple system (the center edge of a windmill $W$ in a block $B\in S_i$ would force the other three edges of $W$ to $B$, similarly as in Theorem~\ref{exW}). Thus $K_n^3$ is decomposed into $m=(1+o(1))n/3$ $W$-free triple systems, showing $r_k(W)\ge (1-o(1))3k$.

The upper bound follows from Theorem~\ref{exW}: in a $k$-coloring of $K_n^3$ with no monochromatic $W$, each color class has at most $\ex(n,W)={n\choose 2}$ edges. Thus ${n\choose 3}/k\le {n\choose 2}$, implying $n\le 3k+2$.  So  by Lemma~\ref{general}(i), $r_k(W)\le 3k+3$.
\qed

We need the following result for tight path.

\begin{proposition}~\label{exP3}
$ex(n,P_3^3)\le {n(n-1)\over 3}$ with equality for  $n\equiv 1,4\mbox{ (mod 12)}$.
\end{proposition}

\noindent {\bf Proof of Proposition.} For a $P_3^3$-free hypergraph $T$ on $n$ vertices and a vertex $v$, the degree $d(v)\le ex(n-1,P_4)\le n-1$. Thus $3|E(T)|=\sum_v d(v)\le n(n-1)$. The statement for equality comes from a $2-(n,4,1)$ design by replacing all blocks by $K_4^3$-s. \qed

\noindent {\bf Proof of Theorem \ref{tightp}.} Observe that the trace of $P_3^3$ at its transversal vertex is $P_4$, the path on four vertices. The upper bound can be obtained in two ways.

Applying and Lemma~\ref{general} (i) with proposition~\ref{exP3}, we have $r_k(P_3^3)\le 2k+3$. On the other hand, we may apply Lemma~\ref{single-transversal} as well: $r_k(P_3^3)\le r_k(P_4)+1\le 2k+3$ (\cite{RA}).

For the lower bound we start with a $3-(n,4,1)$ design $F$ (already used in the proof of Theorem \ref{bow}) and follow the construction in the proof of Theorem \ref{windmill}.  Consider the $6$-uniform hypergraph $H$  with vertex set being the set of pairs  of vertices of $F$ and edges formed by the sets of pairs within the blocks of $F$.  The degree of any vertex in $H$ is $d=(n-2)/2$, the codegree of any pair of vertices is at most one,  so  the conditions for Pippenger-Spencer Theorem are satisfied, giving a decomposition of $H$ into $(1+o(1))d=(1+o(1))n/2$ matchings, $M_i$. Each $M_i$ corresponds to a set $F_i$ of blocks of $F$, intersecting each other in at most one element. Let $T_i$ be the set of triples covered by the blocks of $F_i$. The $T_i$-s provide the required $P_3^3$-free coloring of $K_n^3$ with $(1+o(1))n/2$ colors. \qed
~\\

\section{Small Ramsey numbers}\label{small}
The only known non-trivial classical Ramsey number for triples is $r_2(K_4^{3})=13$, due to McKay and Radziszkowski \cite{MR}.

It is proven in (\cite{RA}  that $13\leq r_3(K_4^{(3)}- e)\leq 16$ and  stated as an easy fact without a proof that  $r_2(K_4^3-e)=7$.
Here we prove this for completeness.

\begin{proposition}  $r_2(K_4^3-e)=7$.
\end{proposition}
\noindent {\bf Proof.}
Consider  the following coloring $C$ of $K_6^3$.
Fix the set of five vertices, $V$,  and let $c$ be the $2$-coloring of $K_5$ on vertex set $V$ with two monochromatic $C_5$'s. Let $v$ be the remaining vertex of $K_6^3$.
 For any triple containing $v$, let $C(\{v,u,w\})=c(uw)$.\\
 For each triple $xyz$, not containing $v$, let $C(\{x,y,z\})$ be the  color different from $c(V-\{x,y,z\})$.
Under coloring $C$, there are two triples of each color in every $4$-set, hence there is no monochromatic $K_4^3-e$.
\qed

The following proposition determines the small undecided cases from Theorem \ref{kite}. A hypergraph is linear if every two edges share at most one vertex.

\begin{proposition} \label{r5kite} $r_4(kite)=5, r_5(kite)=6, r_6(kite)=8$.
\end{proposition}

\noindent {\bf Proof.}
It is obvious that $r_4(kite)>4$. The fact that  $r_4(kite)\le 5$ follows by observing that any $4$-coloring of the edges of $K_5^3$ contains three edges of the same color.

Coloring the triple  $ijk$,  $1\le i<j<k\le 5$ by color $i+j+k \mbox{ (mod 5)}$ gives  $r_5(kite)>5$.
To show that $r_5(kite)\le 6$, we need the result of Cayley \cite{CA}, stating that the maximum number of pairwise disjoint Fano planes in $K_7^3$ is $2$. Suppose $K_6^3$ on vertex set $V$ is $5$-colored so that each color class $i$ is a linear hypergraph $P_i$. Since the average number of edges in a color class is four and no linear hypergraphs on 6 vertices can have more than four edges, it follows that each $P_i$ must be a Pasch configuration. Therefore the pairs uncovered by the triples of $P_i$ form a matching $M_i$ in the complete graph on $V$. The $M_i$-s must form a factorization on $V$ otherwise some pair in $V$ would be covered by at most three $P_i$-s instead of the required four.
These $P_i$-s can be extended by a new vertex to a decomposition of $K_7^3$ into five Fano planes, contradicting Cayley's theorem stated above.

\begin{figure}[ht]
\begin{center}
\includegraphics{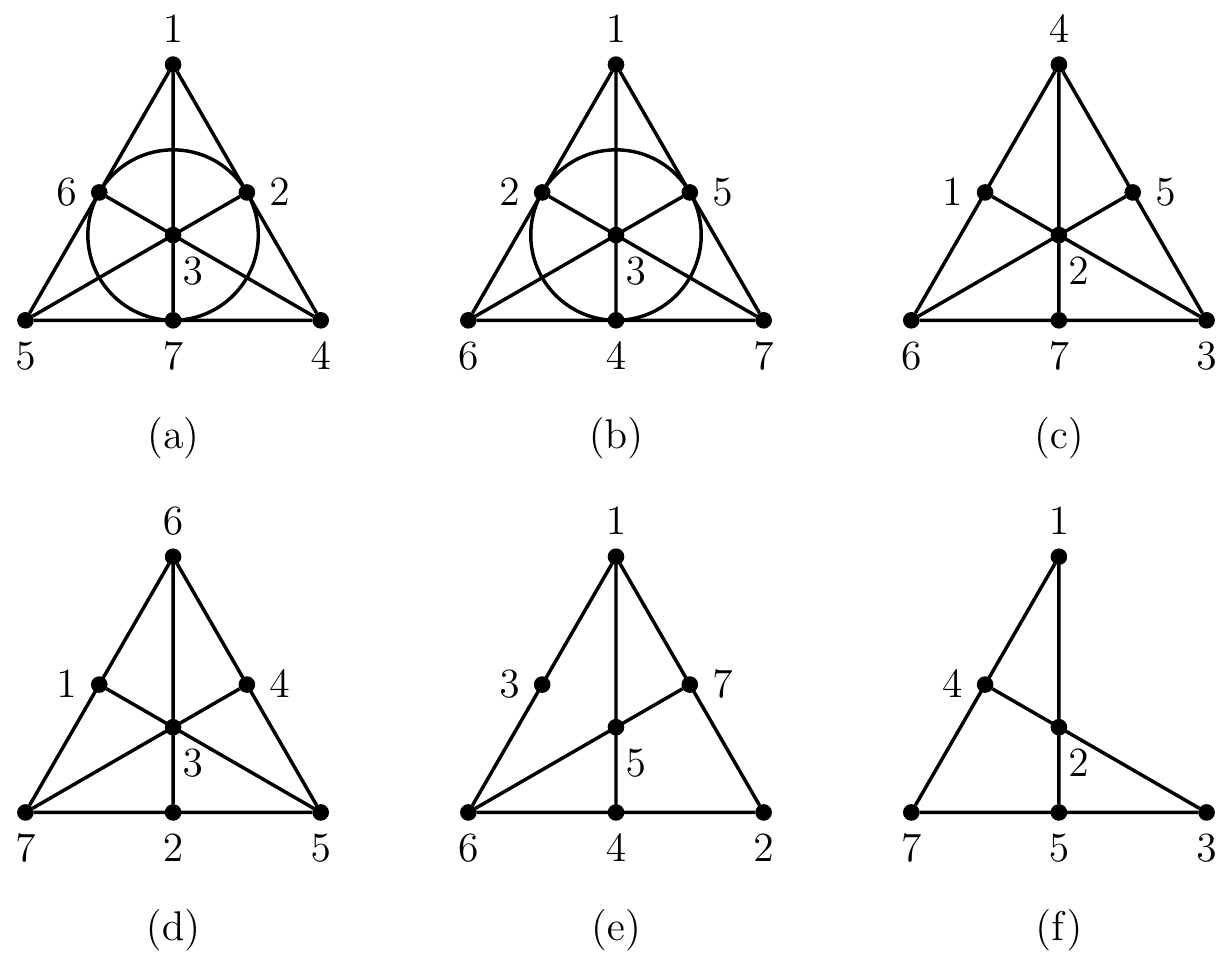}
\end{center}
\caption{Partition of $K_7^3$ into two Fano, two Fano$-e$, Fano$-2e$, Pasch}
\label{fano}
\end{figure}

The upper bound $r_6(kite)\le 8$ is already proved (see the proof of Theorem \ref{kite}). For the lower bound we need a partition of $K_7^3$ into six linear hypergraphs, see Figure~\ref{fano}. Set $V=[7]$ and let $F_1,F_2$ be the two Fano planes generated by shifts of $124,134$ $\pmod{7}$. The next two sets are isomorphic to a Fano plane from which one line is deleted:
$$F_3=\{135,167,236,257,347,456\}, F_4=\{123,146,247,256,345,367\}$$
and $F_6=\{127,136,145,246,567\}$ (Fano plane from which two lines are deleted), $F_7=\{125,147,234,357\}$ (a Pasch configuration).

\qed

\begin{proposition}Set $r_k=r_k(bow)$, then $r_1=r_2=5,r_3=r_4=r_5=6,$ $ r_6=7,r_7=r_8=r_9=r_{10}=9,$
$9\le r_{11}\le r_{12}\le r_{13}\le r_{14}\le 10, r_{15}=11.$
\end{proposition}

\noindent {\bf Proof.}
All but one upper bounds are obtained from Lemma~\ref{general}(i). The exceptional case is when Lemma~\ref{general}(i)  gives $r_5(bow)\le 7$. Here we improve it as follows. Suppose $K_6^3$ is $5$-colored without monochromatic bow. From Lemma \ref{tur} each color class is either a $K_4^3$ (type A) or four triples pairwise intersecting in the same base pair (type B). There are at most three type A colors. The base pairs for different type B colors must be vertex disjoint. Thus there are at least two type A color classes, w.l.o.g. $abcd,cdef$. But then only the base pairs $ae,af,be,bf$ are available for type B colors. Therefore we have two type B and three type A colors, the third is the $K_4^3$ spanned by $abef$. Now there is no base pair available for type B color classes since every pair of vertices is covered by a type A $K_4^3$.

Lower bounds should be exhibited for $r_1,r_3,r_6,r_7,r_{15}$ only. Coloring all triples of $K_4^3$ with the same color, $r_1>4$ follows. Coloring the triples of $\{1,2,3,4\}$ with color $1$, the triples $125,135,235$ with color $2$, the triples $145,245,345$ with color $3$, $r_3>5$ follows. Then $r_6>6$ comes from the following $6$-coloring with color classes
$\binom{\{1,2,3,4\}}{3}$,
 $\binom{\{3,4,5,6\}}{3}$,
 $\binom{\{1,4,5,6\}}{3}-\{4,5,6\}$,
 $\binom{\{2,4,5,6\}}{3}-\{4,5,6\}$,
$\binom{\{1,2,3,5\}}{3}-\{1,2,3\}$,
$\binom{\{1,2,3,6\}}{3}-\{1,2,3\}$.
The $7$-coloring of $K_8^3$ is the $7$ parallel classes of the unique $3-(8,4,1)$ design. Finally, the $15$-coloring of $K_{10}^3$ comes from the unique $3-(10,4,1)$ design whose $30$ blocks can be partitioned into $15$ disjoint pairs.
\qed

\begin{proposition} $r_2(F_5)=6$.
\end{proposition}
\noindent {\bf Proof.}
The lower bound is obvious, color triples of $K_5^3$ containing a fixed vertex with color $1$ and other triples by color $2$. For the upper bound, consider a $2$-colored $K_6^3$ on vertex set $\{1,2,3,4,5,6\}$ and its $2$-colored trace $K=K_5^2$ with respect to vertex $6$. There is a monochromatic, say red odd cycle $C$ in $V(K)-\{6\}$. If $C=1,2,3,1$ then either there is a red triple in $K$ with two vertices on $C$ and
one vertex not in $C$ or all such triples are blue. The former gives a red, the latter a blue $F_5$. If $C=1,2,3,4,5,1$ then either there is a red triple with vertices non-consecutive on $C$ or all the five such triples are blue. Again, the former gives a red, the latter a blue $F_5$.
\qed

\begin{theorem}\label{r3f5} $r_3(F_5)=7$.
\end{theorem}

\noindent {\bf Proof.}
For the  lower bound,  color the triples of $K_6^3$ containing $v$ with color $1$, color uncolored triples containing vertex $w\neq v$ with color $2$ and color all other edges with color $3$.

To prove the upper bound, call a graph $G$ {\it nice} if for every triple
  $T=\{v_1,v_2,v_3\}$ of vertices at least one of the following holds:\\
1. There are two vertex disjoint edges of $G$, such that one of them is in $T$ and the other meets $T$.
2. There is a path of length two in $G$ connecting two vertices of $T$ with midpoint not in $T$.

\begin{observation}\label{trace} If $H$ is an $F_5$-free $3$-uniform hypergraph,  such that the trace of $v$ for a vertex $v$   is   a nice graph,  then all edges of $H$ within $V(G)\cup \{v\}$ contain $v$.
\end{observation}

Indeed, otherwise from the definition of a nice graph we find $F_5$ in $H$. Thus finding a large nice subgraph in a trace one can reduce the number of colors. More generally, a graph is $i$-nice if the property holds for all but at most $i$ triples of vertices.

We need a lemma on $6$-vertex graphs. Since its proof is routine but lengthy, we state it without proof.
\begin{lemma}\label{six} Suppose $G$ has six vertices. If $|E(G)|\ge 9$ then $G$ is nice. If $|E(G)|=8$ then $G$ is $1$-nice, if $|E(G)|=7$ then  $G$ is $2$-nice. If $|E(G)|=6$ then $G$ is $5$-nice, except in one case, when $G$ is $K_{2,3}$ plus an isolated vertex (in this case it is $6$-nice).
\end{lemma}

With these preparations we are ready to prove the upper bound. The majority color, say red in a $3$-colored $K_7^3$, has at least $12$ edges. Some vertex $v$ has red degree at least $6$. Let $G$ be the trace of a red hypergraph at $v$. We get a contradiction from Lemma \ref{six} (and from the fact that we have $12$ edges) except when $G$ has exactly six edges and the trace is $K_{2,3}+w$. This case implies that the red color class has $12$ edges forming  $K_{2,2,3}$, a complete $3$-partite hypergraph with parts of sizes $2$, $2$, and $3$. However, among the $35-12=23$ edges of other colors, one color, say blue, has at least $12$ edges. Repeating the argument for the blue hypergraph, we conclude that the blue hypergraph is also a $K_{2,2,3}$. However, as one can easily check, there is no way to place two edge disjoint $K_{2,2,3}$-s on $7$ vertices.
\qed

\section{Concluding remarks}\label{conclude}

We determined, for $3$-uniform hypergraphs, $r_k$ ranges from $\sqrt{k}$ to double exponential in $k$, and showed a jump in $r_k$ when $H$ changes from $r$-partite to non-$r$-partite. This leads to the following question.

\begin{problem} For which $3$-uniform hypergraphs $F$, is $r_k(F)$ double
exponential? Are there other jumps that the Ramsey function $r_k$ exhibits?
\end{problem}

The ramsey-numbers $r_k(bow),r_k(kite)$ are closely connected to block designs. In case of the kite the only uncertainty is whether $r_k(kite)$ is $k+1$ or $k+2$ when $k\equiv 4\mbox{ (mod 6)}$. This leads to the following problem.

\begin{problem}\label{stspacking}
Suppose $n \equiv 5\mbox{ (mod 6)}$. Is it possible to partition the triples of an $n$-element set into $n-1$ partial triple systems, i.e. into parts so that distinct triples in each part intersect in at most one vertex?  By Theorem \ref{kite}, this is not possible for $n=5$ but perhaps for large enough $n$ (possibly for $n\ge 11$) such partitions exist.
\end{problem}

In case of the bow, the problems related to sharper bounds of $r_k(bow)$ are not purely design theoretic, since color classes can be star components as well. We state just one of those problems.

\begin{problem}\label{k4packing} Suppose $n \equiv 6,10 \mbox{ (mod 12)}$. Is it possible to partition the triples of an $n$-element set into ${n(n-1)\over 2}$ classes so that each class is the union of some disjoint $K_4^3$-s and at most one star component? (Any color class has $n-2$ triples.) For $n=6$  there is no solution.
\end{problem}

Concerning $r_k(K_3-e)$ the most challenging (perhaps difficult) problem is to decrease the upper bound of Theorem \ref{K43-} by one.

\begin{problem}\label{k4packing} $r_k(K_4^3-e)<r_k(K_3)+1$ for every $k\ge 3$?
\end{problem}

A challenging open problem is to improve the estimates of $r_k(P)$ (and/or $ex(n,P)$) where $P$ is the {\em Pasch configuration} with edges $\{abc,bde,cef,adf\}$. (It can be obtained from the Fano plane by deleting a vertex.) Presently only the following is known.

\begin{proposition}\label{pasch}
For positive constants $c, c'$, ~ $c\left(\frac{k}{\log k}\right)^2 < r_k(P) < c'k^4$.
\end{proposition}

\noindent {\bf Proof.} The lower bound is based on the following $P$-free hypergraph, showing that $\ex(n,P)=\Omega(n^{5/2})$, \cite{LPR}.
Take an incidence graph $G$  of a projective plane with $n$ points and $n$ lines. It has  $\Omega(n^{3/2})$ edges. Add $n$ new vertices $x_1,..., x_n$ and add all triples of the form $x_i\cup e$, where $e$ is an edge of $G$.
 The resulting $3$-uniform hypergraph, call it $H$, has $3n$ vertices and $\Omega(n^{5/2})$ edges.

Notice that the edge-density of $H$ is $d(H)=cn^{-1/2}$ for some constant $c>0$.
From Lemma~\ref{general}(ii) we see that there is a coloring of $K_n^3$ with $(c'n^{1/2}\log n)$ colors and no monochromatic $P$.
Thus $r_k(P)>n$ with $k=c'n^{1/2}\log n$.  Expressing $n$ in terms of $k$ gives  the desired lower bound.

The upper bound follows from  Lemma~\ref{general}(i) and  the fact that $\ex(n,P)=O(n^{11/4})$ \cite{LPR}. This is based on the claim that $\ex(n,K(2,2,2))=O(n^{11/4})$ proved by Erd\H os \cite{E}, where $K(2,2,2)$ is the complete $3$-partite $3$-uniform hypergraph with two vertices in each part.
 \qed

\section{Acknowledgments}

Thanks to Zoli F\"uredi  and Roman Glebov for conversations on the subject of this paper, and Stefan Walzer for improving the lower bound in Theorem \ref{thm-complete}. The authors would also like to thank the organizers of the  2012 Midsummer Combinatorial Workshop that took place at  Charles University.

\end{document}